\documentclass[12pt,a4paper]{article}
\usepackage[utf8]{inputenc}
\usepackage[english]{babel}
\usepackage{amsfonts,amssymb,amsthm,setspace}
\usepackage{graphicx,anysize}
\usepackage{latexsym,amsmath,authblk,hyperref}
\usepackage{array,xcolor,multirow,multicol,lscape,setspace}
\usepackage[margin=1in]{geometry} 
\usepackage[justification=centering]{caption}

\theoremstyle{definition}

\begin{document}
\title{The extended xgamma distribution}
\begin{scriptsize}
\author{Mahendra Saha $^{a}$, Abhimanyu Singh Yadav $^{a}$\footnote{Corresponding author. e-mail: abhmanyu@curaj.ac.in}, Arvind Pandey $^{a}$, Shivanshi Shukla $^{a}$ and Sudhansu S Maiti $^{b}$\\
\small $^{a}$Department of Statistics, Central University of Rajasthan, Rajasthan, India\\
\small $^{b}$ Department of Statistics, Visva-Bharati University, Santiniketan, India}
\end{scriptsize}
\date{}
\maketitle
\begin{abstract}
This article aims to introduced a new distribution named as extended xgamma (EXg) distribution. This generalization is derived from xgamma distribution (Xg), a special finite mixture of exponential and gamma distributions [see, Sen et al. ($2016$)]. Some important statistical properties, viz., survival characteristics, moments, mean deviation and random number generation have been derived. Further, maximum likelihood estimation for the estimation of the unknown parameters have also been discussed for the complete sample. The application of the proposed model has been illustrated through a real data set and observed that the proposed model might be taken as an better alternative to some well known lifetime distributions.  
\end{abstract}

{ \bf Keywords:} Exponential distribution, gamma distribution, xgamma distribution, moments, maximum likelihood estimation.

\section{Introduction}
In reliability analysis, the lifetime of any electronic device or items is varying in nature. Hence, it seems to be logical to model the lifetime data with a specific probability distribution. The exponential distribution and it's different generalizations, e.g., Weibull, gamma, exponentiated exponential etc. have been often used to model the data with constant, monotone hazard rate functions. Also, the finite mixtures of two or more probability distributions are also the better alternative to analyze any life time data, such as Lindley [see, Lineley ($1958$)], generalized Lindley [see, Nadarajah et al.($2011$)]. In the same era of generalization of statistical distributions, the one parameter family of distributions, namely, xgamma (Xg) distribution  is one of them which is a special finite mixture of exponential and gamma distributions, proposed by Sen et al. ($2016$). The probability density function (PDF) of Xg distribution is given as
\begin{eqnarray*}
f(x;\theta)&=&\frac{\theta^2}{(1+\theta)}\left(1+\frac{\theta}{2}x^2\right)e^{-\theta x}~~;x>0,~\theta>0
\end{eqnarray*}
They have discussed various mathematical properties, viz., moments, reliability characteristics and stochastic ordering etc. They have also discussed the estimation of the parameter and shown the superiority of Xg distribution over exponential distribution. Yadav et al. ($2018b$), studied the Bayesian estimation of the parameter and the reliability characteristics of Xg distribution using Type-II hybrid censored data. In most of the situations, finite mixture distributions arising from the standard distributions play a better role in modelling lifetime phenomena as compared to the standard distributions. Recently, Yadav et al. ($2018a$) introduced the inverted version of Xg distribution which possesses the upside-down bathtub-shaped hazard function. The Xg distribution did not provide enough flexibility for analyzing different types of lifetime data as it is of one parameter. It will be useful to consider further alternatives of Xg distribution to increase the flexibility for modelling purposes.\\

In this article, we propose a three parameter family of distribution which generalizes the Xg distribution, named as the extended xgamma (EXg) distribution and hence the name proposed. The procedure used is based on certain finite mixtures of exponential and gamma distributions. The shape parameter provides more flexibility for describing different types of data allowing hazard rate modelling. Moreover, we also derived some statistical characteristics such as, survival and hazard rate functions, moments, mean deviation etc. The unknown parameters of the model are estimated via method of maximum likelihood estimation (MLE). Besides, the compatibility of the proposed model has been shown based on a real data set and is observed that EXg distribution is best as compared to its particular models. To the best of our knowledge, no such generalization from Xg distribution has been used thus so for; hence, this study is targeted to mold the gap through this present study. 
\section{Definition}
Let,
\begin{eqnarray}\label{eq1}
f(x;\alpha,\theta)&=&\frac{\theta^\alpha}{\Gamma(\alpha)}e^{-\theta x}x^{\alpha-1}~~;x>0,~\theta>0
\end{eqnarray}
be the probability density function (PDF) of the gamma distribution with shape parameter $\alpha$ and the scale parameter $\theta$ and is denoted by gamma($\alpha$, $\theta$). Let $V_1$ and $V_2$ are two random variables which are distributed according to gamma($\alpha$, $\theta$) and gamma($\alpha+2$, $\theta$) respectively. Suppose for $\beta~(\ge 0)$, consider the random variable $X=V_1$ with probability $\frac{\theta}{\theta+\beta}$ and $X=V_2$ with probability $\frac{\beta}{\theta+\beta}$. Then, it is easy to verify that the PDF of $X$ is 
\begin{eqnarray}\label{eq2}
f(x;\alpha,\theta,\beta)&=&\frac{\theta^{\alpha+1}}{(\theta+\beta)\Gamma(\alpha+2)}\left(\alpha^2+\alpha+\theta\beta x^2\right)e^{-\theta x}x^{\alpha-1}~~;x>0,~\alpha,~\theta,~\gamma>0
\end{eqnarray}
As a particular case, when $\alpha=\beta=1$, the distribution, given in Equation ($\ref{eq2}$), contains the Xg distribution. The Equation ($\ref{eq2}$) reduces to the PDF of the gamma distribution (GD) with parameters $\alpha$ and $\theta$ when $\beta=0$. Again, the case, where $\alpha=1$ and $\beta=0$, Equation ($\ref{eq2}$) coincides with the PDF of the exponential distribution (ED). Hence, we say that the random variable $X$ has a EXg distribution  with parameters $\alpha$, $\theta$ and $\beta$, donated as EXg($\alpha,\theta,\beta$). To find the cumulative distribution function (CDF) of $X$, where, $X$ is following EXg($\alpha,\theta,\beta$), we have to calculate
\begin{eqnarray}\label{eq3}
F(x;\alpha,\theta,\beta)&=&P(X \le x)\nonumber
\\&=&\frac{\theta^{\alpha+1}}{(\theta+\beta)\Gamma(\alpha+2)}\left[\int\limits_{0}^{x}\left(\alpha^2+\alpha+\theta\beta u^2\right)e^{-\theta u}u^{\alpha-1} du \right]\nonumber
\\&=&\frac{(1+\alpha)}{\Gamma(\alpha-1)}\gamma_l(\alpha,\theta x)+\frac{\beta \theta}{\Gamma(\alpha)}\gamma_l(\alpha+2,\theta x)~;~~~\alpha>1
\end{eqnarray} 
where, $\gamma_l(\alpha,\theta x)$ is the lower incomplete gamma function.

\begin{figure}
	\includegraphics[width=3in,height=4in]{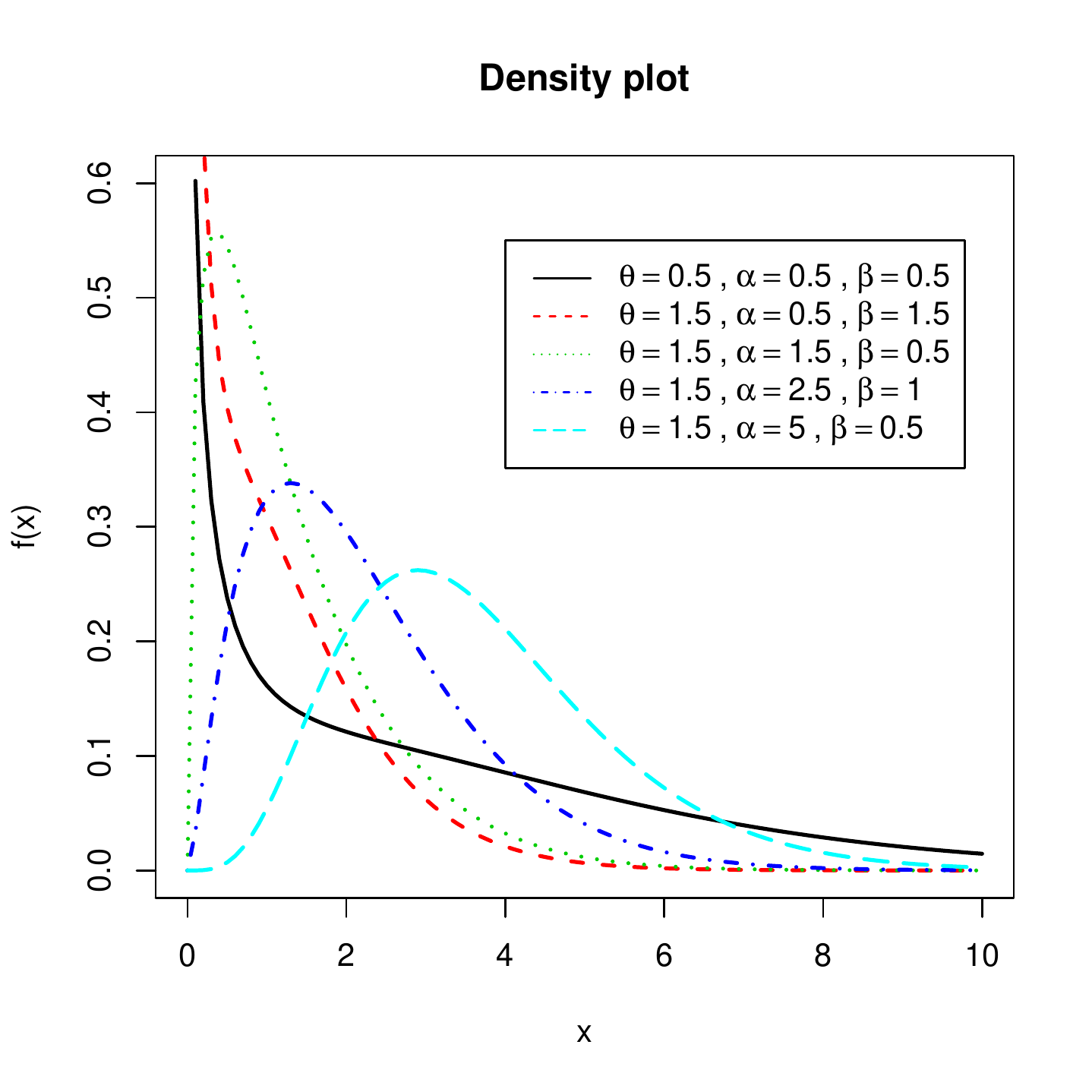}
		\includegraphics[width=3in,height=4in]{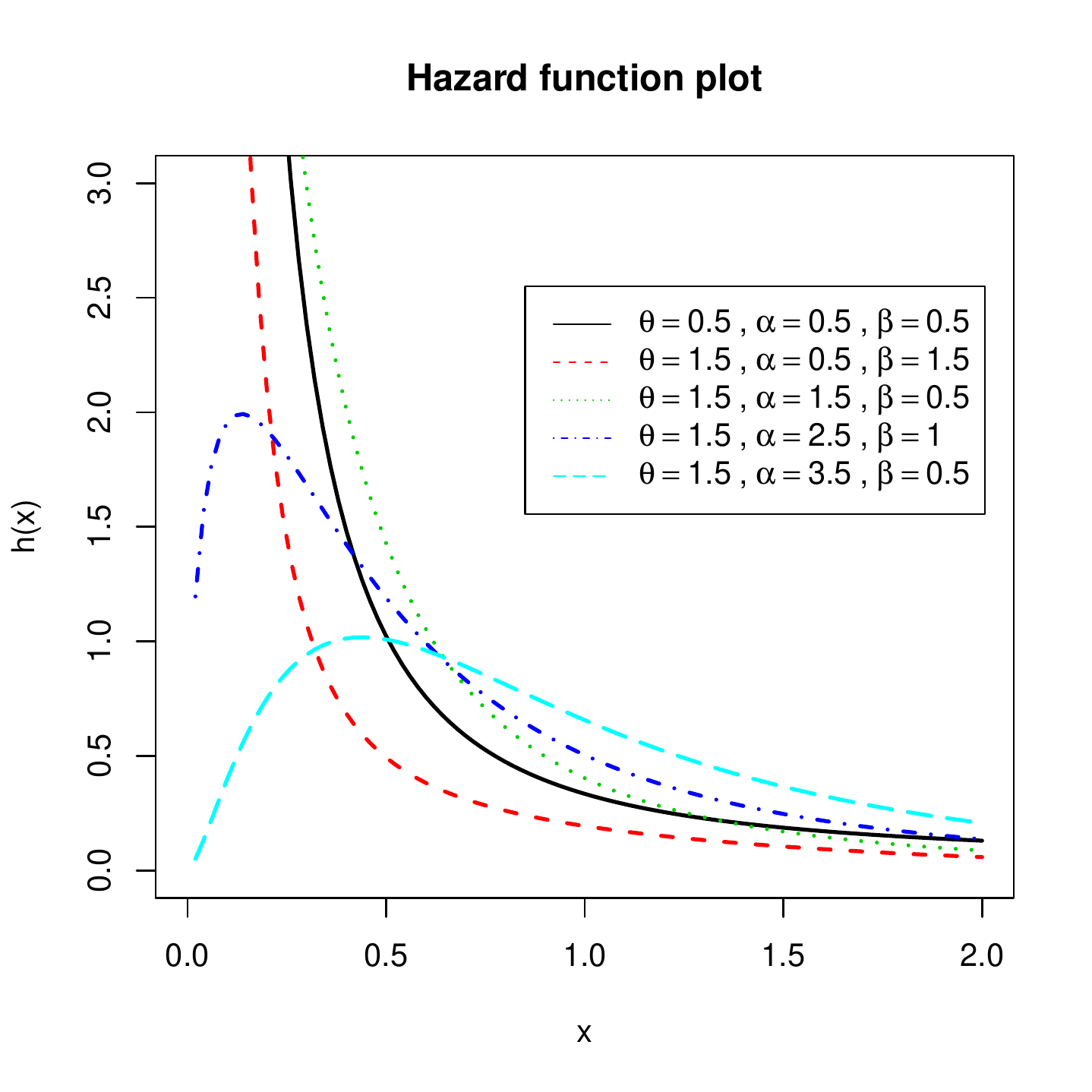}
		\caption{Density and hazard function plot.}
\end{figure}
Shape of the distribution is traced graphically for the PDF ($2$) and it is observed that the proposed model is highly positively skewed and uni-model, [see Figure 1]. The shape is highly depended over the values of $\alpha$, $\beta$ and for higher $\alpha$ shape tends to symmetric distribution.
%
\subsection{Survival and hazard rate functions}
The survival and hazard function are the two most important lifetime characteristics and frequently used to describe the inherent properties of several survival data. The survival function is the probability of an individual or any lifetime system which survive beyond the mission time $t_0$.  Let $T$ be the lifetime random variable  follows PDF (2). Let survival function and hazard function are denoted by $S(t_o)$ and $h(t_o)$, then the survival and hazard functions for EXg distribution are given by the following expressions;

\begin{equation}
S(t_o)=\frac{(1+\alpha)}{\Gamma(\alpha-1)}\gamma_u(\alpha,\theta x)+\frac{\beta \theta}{\Gamma(\alpha)}\gamma_u(\alpha+2,\theta x)
\end{equation}

\begin{equation}
h(t_o)=\dfrac{\frac{\theta^{\alpha+1}}{(\theta+\beta)\Gamma(\alpha+2)}\left(\alpha^2+\alpha+\theta\beta x^2\right)e^{-\theta x}x^{\alpha-1}}{\frac{(1+\alpha)}{\Gamma(\alpha-1)}\gamma_u(\alpha,\theta x)+\frac{\beta \theta}{\Gamma(\alpha)}\gamma_u(\alpha+2,\theta x)}
\end{equation}
respectively and $\gamma_l(\alpha,\theta x)$ represents the lower incomplete gamma function. The shape of the hazard function for different combination of the model parameters has been presented in Figure 1 and is observed that the proposed distribution possess increasing, decreasing and hump type hazard rate. 

\subsection{Random variate generation}
The PDF of the EXg distribution can be written in terms of finite mixture of the PDF of gamma distribution as
\begin{eqnarray*}
f(x;\alpha,\theta,\beta)&=&\frac{\theta}{\theta+\beta}f(x;\alpha,\theta)+\frac{\beta}{\theta+\beta} f(x;\alpha+2,\theta)
\end{eqnarray*} 
To generate random number from EXg($\alpha,\theta,\beta$), the following steps may be used.\\
\begin{enumerate}
\item Specified the values of $\alpha$, $\theta$, $\beta$ and $n$.
\item Generate $U_i$ from $uniform(0,1)$ distribution $(i=1,2,...,n)$.\\
\item Generate $V_i$ from $gamma(\alpha,\theta)$ distribution $(i=1,2,...,n)$.\\
\item Generate $W_i$ from $gamma(\alpha+2,\theta)$ distribution $(i=1,2,...,n)$.\\
\item  If $U_i\leq\frac{\theta}{\theta+\beta}$, set $X_i=V_i$, otherwise set $X_i=W_i$.\\
\end{enumerate}
If we take $\alpha=1$, then we get the random variates from Xg distribution.

\section{Some important statistical properties}
In this section, we have studied some important statistical properties EXg distribution.
\subsection{Moments}
Here, we have obtained the expression of the $r$th order moment about origin, i.e., $r$-th order raw moment for the EXg distribution.
\begin{eqnarray}
E(X^r)&=& \int\limits_{x} x^r f(x)dx\nonumber
\\&=&\frac{\theta^{\alpha+1}}{(\theta+\beta)\Gamma(\alpha+2)} \int\limits_{0}^{\infty} \left(\alpha^2+\alpha+\theta\beta x^2\right)e^{-\theta x}x^{\alpha+r-1} dx\nonumber
\\&=&\frac{\Gamma(\alpha+r)}{\theta^{r-1}(\theta+\beta)\Gamma(\alpha+2)}\left[\alpha^2+\frac{\beta}{\theta}(\alpha+r)(\alpha+r+1)\right]~;~r=1,2,3,...
\end{eqnarray}
Now, in particular, for $r=1$, we just get the expression for first order raw moment, which is nothing but mean of $X$, where, $X$ is following EXGD, i.e.,
\begin{eqnarray}
E(X)&=&\frac{1}{(\theta+\beta)(\alpha+1)}\left[\alpha^2+\frac{\beta}{\theta}(\alpha^2+3\alpha+2)\right]
\end{eqnarray}
Also, for $r=2,~3~\mbox{and}~4$, we get the second, third and fourth order raw moments respectively, given as
\begin{eqnarray}
E(X^2)&=&\frac{1}{\theta(\theta+\beta)}\left[\alpha^2+\frac{\beta}{\theta}(\alpha+2)(\alpha+3)\right]
\end{eqnarray}
\begin{eqnarray}
E(X^3)&=&\frac{(\alpha+2)}{\theta^2(\theta+\beta)}\left[\alpha^2+\frac{\beta}{\theta}(\alpha+3)(\alpha+4)\right]
\end{eqnarray}
\begin{eqnarray}
E(X^4)&=&\frac{(\alpha+2)(\alpha+3)}{\theta^3(\theta+\beta)}\left[\alpha^2+\frac{\beta}{\theta}(\alpha+4)(\alpha+5)\right]
\end{eqnarray}
Again, by using the relation between the raw moments and central moments, we can find the expression for successive central moments for EXg distribution. Hence, second order central moment, i.e., variance is given as, where $X$ is following EXg distribution
\begin{eqnarray}
\mbox{Variance}(X)&=&E(X^2)-E^2(X)\nonumber
\\&=&\frac{1}{\theta(\theta+\beta)}\left[\alpha^2+\frac{\beta}{\theta}(\alpha+2)(\alpha+3)\right]-\nonumber
\\&&\left[ \frac{1}{(\theta+\beta)(\alpha+1)}\left\{\alpha^2+\frac{\beta}{\theta}(\alpha^2+3\alpha+2)\right\} \right]^2
\end{eqnarray} 
third and fourth order central moments can also be obtained in the same lines. We can also obtained the coefficient of variation (CV) from the above expressions, given as
\begin{eqnarray}
\mbox{CV}&=&\frac{\sqrt{\mbox{Variance}(X)}}{E(X)}\nonumber
\\&=&\frac{\alpha\left[\left\{\frac{1}{\theta}+\frac{\beta}{\theta^2}(1+\frac{2}{\alpha})(1+\frac{3}{\alpha})\right\}-\frac{\alpha^2}{(\theta+\beta)(\alpha!)^2}\left\{1+\frac{\beta}{\theta}(1+\frac{2}{\alpha})(1+\frac{1}{\alpha})\right\}\right]}{\sqrt{\theta+\beta}(\alpha+1)\left[\alpha^2+\frac{\beta}{\theta}(\alpha+1)(\alpha+2)\right]}
\end{eqnarray}
\subsection{Mean deviation}
The mean deviation about mean of random variable $X$, having PDF, given in Equation ($\ref{eq2}$) is obtained as
\begin{eqnarray}
MD_{\mu}&=& E|X-\mu|\nonumber
\\&=&\int\limits_{0}^{\infty}|(x-\mu)| f(x) dx \nonumber
\\&=&\int\limits_{o}^{\mu}(\mu-x)f(x)dx+\int\limits_{\mu}^{\infty}(x-\mu)f(x)dx \nonumber
\\&=&2\mu F(\mu)-\mu+2\int\limits_{\mu}^{\infty}x f(x)dx \nonumber
\\&=&\frac{2}{(\theta+\beta)\Gamma(\alpha+2)}\sum\limits_{r}\frac{(-\theta)^{r+\alpha+1}}{r!}\mu^{r+\alpha+1}\left[\frac{\alpha(\alpha+1)}{(r+\alpha)(r+\alpha+1)}+\frac{\mu^2\beta\theta}{(r+\alpha+2)(r+\alpha+3)}\right]
\end{eqnarray}
where, $\mu=E(X)$ and $F(\mu)$ stands for CDF of X upto point $\mu$. 
\section{Estimation of the parameters}
Here, we have estimated survival and hazard estimates of proposed distribution EXg distribution by using the maximum likelihood method of estimates using the invariance properties of that. Suppose $X_1,X_2,\cdots,X_n$ be random sample drawn from EXg distribution and $x_1,x_2,\cdots,x_n$ be the observed values from the original sample. Then, the likelihood function is given as 
\begin{eqnarray}\label{eq12}
L(\alpha,\theta,\beta\mid x)&=&\prod\limits_{i=1}^{n}\left[\frac{\theta^{\alpha+1}}{(\theta+\beta)\Gamma(\alpha+2)}\left(\alpha^2+\alpha+\theta\beta x_i^2\right)e^{-\theta x_i}x_i^{\alpha-1}\right]
\end{eqnarray}
Taking logarithm on both the sides of the above Equation ($\ref{eq3}$), we just get the log-likelihood function and is given as
\begin{eqnarray}\label{eq13}
\log L(\alpha,\theta,\beta\mid x)&=&n(\alpha+1)\log(\theta)-n \log{(\theta+\beta)}-n \log{\Gamma(\alpha+2)}+\nonumber
\\&& \sum\limits_{i=1}^{n} \log\left(\alpha^2+\alpha+\theta\beta x_i^2\right)+\sum\limits_{i=1}^{n}(\alpha-1) \log{x_i}-\theta\sum\limits_{i=1}^{n}x_i
\end{eqnarray}
Now, we know that partial derivatives of the log-likelihood functions with respect to $\alpha$, $\theta$ and $\beta$ and equating to zero, yield the MLEs $\hat{\alpha}$, $\hat{\theta}$ and $\hat{\beta}$, i.e.,
\begin{eqnarray}
\frac{\partial \log L(\alpha,\theta,\beta\mid x)}{\partial\alpha}=0
\end{eqnarray}
\begin{eqnarray}
\frac{\partial \log L(\alpha,\theta,\beta\mid x)}{\partial\theta}=0
\end{eqnarray}
and
\begin{eqnarray}
\frac{\partial \log L(\alpha,\theta,\beta\mid x)}{\partial\beta}=0
\end{eqnarray}
Since, the above non-linear equation can not be solved analytically; thus any iterative procedure has been used to obtain the estimates of the parameters. 


\section{Real data analysis}
In this Section, we have considered one data set which represents the strength of $1.5$ cm glass fibers measured at the National Physical Laboratory, England. Unfortunately, the units of measurements are not given in the paper and are taken from Smith and Naylor ($1987$). The fitting of the proposed model has been compared with the following lifetime distributions.
\begin{itemize}
	\item Exponential distribution (ED)
    \item Gamma distribution (GD)
    \item Xgamma (Xg) distribution
\end{itemize} 
The fitting/compatibility has been performed using different model selection tools, namely, negative of log-likelihood, one sample K-S statistic and corresponding $p$ value. The model with minimum negative of log-likelihood,  K-S statistic and maximum of $p$ value is treated as best model. The obtained measures are reported in the following table which indicates that the EXg is best choices among the considered probability distributions. hence EXg may be chosen as an alternative model. All the computation related to the fitting of the data has been performed using $R-$ software, see Ikha and Gentelman (1996).
\begin{table}[htbp]
	\centering
	\caption{Values of the estimate of the parameters and model selection tools.}
	\begin{tabular}{ccccccc}\\
		\hline
		Model & Estimates   &-2 Log-Likelihood    & K-S statistic  \\
		\hline
		ED    & $\hat{\theta}=0.6636$ & $177.6600$ & $0.5640$  \\
		GD    & [$\hat{\theta}=11.5711$, $\hat{\alpha}=17.4355$] & $47.9030$ & $0.2164$   \\
		Xg   & $\hat{\theta}$=1.3376 &171.6012 & 0.4037   \\
     EXg  & [$\hat{\theta}=0.0156$, $\hat{\alpha}=17.5700$, $\hat{\beta}=11.6800$] & $43.5070$& 0.1951 \\
	\hline
	\end{tabular}%
	\label{tab:addlabel}%
\end{table}%
\section{Concluding remarks}
In this article, a new three parameter model as an extension of Xg distribution has been proposed and studied. Some statistical properties such as survival characteristics, moments, mean deviation about mean etc. have been derived and discussed. The algorithm of random number generation of this proposed model  has also been given. Further, MLEs of the parameters are also obtained using iterative procedure. Lastly, a real data set has been used to demonstrate the practical applicability of the proposed model and observed that EXg model provides better fit as compared to its particular models. Hence, the proposed model might be taken as an alternative model analyze several reliability/ survival data. 

\end{document}